\documentclass[psamsfonts,reqno]{amsart}
\usepackage{amssymb,eucal,graphics,xy}
\usepackage{epstopdf}
\xyoption{all}

\newcommand{\cml}{complemented modular lattice}
\newcommand{\scml}{sectionally complemented modular lattice}
\newcommand{\Ban}{Ba\-na\-schew\-ski}
\newcommand{\msd}{meet-sem\-i\-dis\-trib\-u\-tive}
\newcommand{\msdy}{meet-sem\-i\-dis\-trib\-u\-tiv\-i\-ty}

\numberwithin{equation}{section}

\newcommand{\pup}[1]{\textup{(}#1\textup{)}}

\theoremstyle{plain}

\newtheorem{lemma}{Lemma}[section]
\newtheorem{theorem}[lemma]{Theorem}
\newtheorem{proposition}[lemma]{Proposition}
\newtheorem{corollary}[lemma]{Corollary}

\newtheorem{claim}{Claim}

\newtheorem*{sclaim}{Claim}

\newtheorem*{stat}{\name}
\newcommand{\name}{testing}

\theoremstyle{definition}
\newtheorem{definition}[lemma]{Definition}
\newtheorem{example}[lemma]{Example}
\newtheorem{problem}{Problem}

\theoremstyle{remark}
\newtheorem{remark}[lemma]{Remark}

\newtheorem*{note}{Note}

\newenvironment{all}[1]{\renewcommand{\name}{#1}\begin{stat}}
                        {\end{stat}}

\newcommand{\qedc}{{\qed}~{\rm Claim~{\theclaim}.}}
\newcommand{\qedsc}{{\qed}~{\rm Claim.}}

\newenvironment{cproof}
{\begin{proof}[Proof of Claim.]}
{\qedc\renewcommand{\qed}{}\end{proof}}

\newenvironment{scproof}
{\begin{proof}[Proof of Claim.]}
{\qedsc\renewcommand{\qed}{}\end{proof}}

\newcommand{\cB}{\mathcal{B}}

\newcommand{\cF}{\mathcal{F}}

\newcommand{\bdu}{{\boldsymbol{u}}}
\newcommand{\bdv}{{\boldsymbol{v}}}

\newcommand{\id}{\mathrm{id}}
\newcommand{\dnw}{\mathbin{\downarrow}}
\newcommand{\upw}{\mathbin{\uparrow}}
\newcommand{\res}{\mathbin{\restriction}}

\newcommand{\jz}{\ensuremath{(\vee,0)}}

\newcommand{\jh}{join-ho\-mo\-mor\-phism}

\newcommand{\eps}{\varepsilon}

\newcommand{\cm}{commutative monoid}

\newcommand{\set}[1]{\{{#1}\}}
\newcommand{\setm}[2]{\set{{#1}\mid{#2}}}

\newcommand{\famm}[2]{({#1}\mid{#2})}

\newcommand{\into}{\hookrightarrow}
\newcommand{\onto}{\twoheadrightarrow}
\newcommand{\ol}[1]{\overline{#1}}

\DeclareMathOperator{\Id}{Id}
\DeclareMathOperator{\Dim}{Dim}

\newcommand{\LL}{\mathbb{L}}
\DeclareMathOperator{\Idemp}{Idemp}
\DeclareMathOperator{\At}{At}

\DeclareMathOperator{\im}{im}

\newcommand{\utr}{\trianglelefteq}

\author[F.~Wehrung]{Friedrich Wehrung}
\address{LMNO, CNRS UMR 6139\\
D\'epartement de Math\'ematiques, BP 5186\\
Universit\'e de Caen, Campus 2\\
14032 Caen cedex\\
France}
\email{wehrung@math.unicaen.fr}
\urladdr{http://www.math.unicaen.fr/\~{}wehrung}

\date{\today}

\begin{document}

\title[Banaschewski functions and coordinatization]{Coordinatization of lattices by regular rings without unit and Banaschewski functions}

\begin{abstract}
A \emph{\Ban\ function} on a bounded lattice~$L$ is an antitone self-map of~$L$ that picks a complement for each element of~$L$. We prove a set of results that include the following:
\begin{itemize}
\item Every \emph{countable} \cml\ has a \Ban\ function with Boolean range, the latter being unique up to isomorphism.

\item Every (not necessarily unital) \emph{countable} von~Neumann regular ring~$R$ has a map~$\eps$ from~$R$ to the idempotents of~$R$ such that $xR=\eps(x)R$ and $\eps(xy)=\eps(x)\eps(xy)\eps(x)$ for all $x,y\in R$.

\item Every \scml\ with a \emph{\Ban\ trace} (a weakening of the notion of a \Ban\ function) embeds, as a neutral ideal and within the same quasivariety, into some \cml. This applies, in particular, to any \scml\ with a countable cofinal subset.
\end{itemize}
A \scml~$L$ is \emph{coordinatizable}, if it is isomorphic to the lattice~$\LL(R)$ of all principal right ideals of a von~Neumann regular (not necessarily unital) ring~$R$. We say that~$L$ has a \emph{large $4$-frame}, if it has a homogeneous sequence $(a_0,a_1,a_2,a_3)$ such that the neutral ideal generated by~$a_0$ is~$L$. J\'onsson proved in 1962 that if~$L$ has a countable cofinal sequence and a large $4$-frame, then it is coordinatizable. We prove that \emph{A \scml\ with a large $4$-frame is coordinatizable if{f} it has a \Ban\ trace}.
\end{abstract}

\maketitle

\section{Introduction}\label{S:Intro}
Bernhard Banaschewski proved in~\cite{Bana} that on every vector space~$V$, over an arbitrary division ring, there exists an \emph{order-reversing} (we say \emph{antitone}) map that sends any subspace~$X$ of~$V$ to a complement of~$X$ in~$V$. Such a function was used in~\cite{Bana} for a simple proof of Hahn's Embedding Theorem that states that every totally ordered abelian group embeds into a generalized lexicographic power of the reals.

By analogy with \Ban's result, we define a \emph{\Ban\ function} on a bounded lattice~$L$ as an antitone self-map of~$L$ that picks a complement for each element of~$L$ (Definition~\ref{D:BanLatt}).
Hence \Ban's above-mentioned result from~\cite{Bana} states that the subspace lattice of every vector space has a \Ban\ function. This result is extended to all \emph{geometric} (not necessarily modular) lattices in Saarim\"aki and Sorjonen~\cite{SaSo}.

We prove in Theorem~\ref{T:BaCML} that \emph{Every countable \cml\ has a \Ban\ function with Boolean range}. We also prove (Corollary~\ref{C:DeltaBVm}) that such a Boolean range is uniquely determined up to isomorphism. In a subsequent paper~\cite{BanCoord2}, we shall prove that the countability assumption is needed.

Then we extend the notion of a \Ban\ function to \emph{non-unital} lattices, thus giving the notion of a \emph{\Ban\ measure} (Definition~\ref{D:BanMeas}) and the more general concept of a \emph{\Ban\ trace} (Definition~\ref{D:BanTail})---first allowing the domain to be a cofinal subset and then replacing the function by an indexed family of elements. It follows from~\cite[Lemma~5.2]{BanCoord2} that every \Ban\ measure on a cofinal subset  is a \Ban\ trace. \Ban\ measures are proved to exist on any countable \scml\ (Corollary~\ref{C:BanTailId2}), and every \scml\ with a \Ban\ trace embeds, as a neutral ideal and within the same quasivariety, into some \cml\ (Theorem~\ref{T:BanTailId}). In particular (Corollary~\ref{C:BanTailId1}),
\begin{quote}\em
Every \scml\ with a countable cofinal subset embeds, as a neutral ideal and within the same quasivariety, into some \cml.
\end{quote}

We finally relate \Ban\ functions to the problem of \emph{von Neumann coordinatization}. We recall what the latter is about. A ring (associative, not necessarily unital) $R$ is \emph{von~Neumann regular}, if for each $x\in R$ there exists $y\in R$ such that $xyx=x$ (cf. Fryer and Halperin~\cite{FrHa54}, Goodearl~\cite{Good91}). The set~$\LL(R)$ of all principal right ideals of a (not necessarily unital) von~Neumann regular ring~$R$, that is,
 \[
 \LL(R):=\setm{xR}{x\in R}=\setm{xR}{x\in R\text{ idempotent}}\,.
 \]
ordered by inclusion, is a sublattice of the lattice of all ideals of~$L$; hence it satisfies the \emph{modular law},
 \[
 X\supseteq Z\ \Longrightarrow\ X\cap(Y+Z)=(X\cap Y)+Z\,.
 \]
(Here~$+$ denotes the addition of ideals.)
Moreover, $\LL(R)$ is \emph{sectionally complemented} (cf. \cite[Section~3.2]{FrHa54}), that is, for all principal right ideals~$X$ and~$Y$ such that $X\subseteq Y$, there exists a principal right ideal~$Z$ such that $X\oplus Z=Y$. A lattice is \emph{coordinatizable}, if it is isomorphic to~$\LL(R)$ for some von~Neumann regular ring~$R$; then we say that~$R$ \emph{coordinatizes}~$L$. In particular, every coordinatizable lattice is sectionally complemented modular. One of the weakest known sufficient conditions, for a \scml, to be coordinatizable, is given by a result obtained by Bjarni J\'onsson in 1960, see~\cite{Jons60}:

\begin{all}{J\'onsson's Coordinatization Theorem}
Every \cml\ $L$ that admits a large $4$-frame, or which is \emph{Arguesian} and that admits a large $3$-frame, is coordinatizable.
\end{all}

We refer to Section~\ref{S:Basic} for the definition of a large $n$-frame. J\'onsson's result extends von Neumann's classical Coordinatization Theorem; his proof has been recently substantially simplified by Christian Herrmann~\cite{Herr}. On another track, the author proved that there is no first-order axiomatization for the class of all coordinatizable lattices with unit~\cite{CXCoord}.

We introduce a ring-theoretical analogue of \Ban\ functions (Definition~\ref{D:BanRing}), and we prove that a unital von~Neumann regular ring~$R$ has a \Ban\ function if{f} the lattice~$\LL(R)$ has a \Ban\ function (Lemma~\ref{L:BanLattRing}). Interestingly, the definition of a \Ban\ function for a ring does not involve the unit; this makes it possible to prove the following result (cf. Corollary~\ref{C:BanRegRings}):

\begin{quote}\em
For every countable \pup{not necessarily unital} von~Neumann regular ring~$R$, there exists a map~$\eps$ from~$R$ to the idempotents of~$R$ such that $xR=\eps(x)R$ and $\eps(xy)=\eps(x)\eps(xy)\eps(x)$ for all $x,y\in R$.
\end{quote}

Finally, we relate coordinatizability of a lattice~$L$ and existence of \Ban\ traces on~$L$. Our main result in that direction is that \emph{A \scml\ that admits a large $4$-frame, or which is Arguesian and that admits a large $3$-frame, is coordinatizable if{f} it has a \Ban\ trace} (Theorem~\ref{T:CharCoord4Fr}).

\section{Basic concepts}\label{S:Basic}

By ``countable'' we will always mean ``at most countable''. We shall denote by~$\omega$ the set of all non-negative integers.

Let~$P$ be a partially ordered set.
We denote by~$0_P$ (resp., $1_P$) the least element (resp. largest element) of~$P$ when they exist, also called \emph{zero} (resp., \emph{unit}) of~$P$, and we simply write~$0$ (resp., $1$) in case~$P$ is understood. Furthermore, we set~$P^-:=P\setminus\set{0_P}$. We set
\begin{align*}
 U\dnw X&:=\setm{u\in U}{(\exists x\in X)(u\leq x)}\,,\\
 U\upw X&:=\setm{u\in U}{(\exists x\in X)(u\geq x)}\,,
 \end{align*}
for any subsets~$U$ and~$X$ of~$P$, and we set
$U\dnw x:=U\dnw\set{x}$, $U\upw x:=U\upw\set{x}$, for any~$x\in P$. We say that~$U$ is a \emph{lower subset} (resp., \emph{upper subset}) of~$P$, if $U=P\dnw U$ (resp., $U=P\upw U$). We say that~$P$ is \emph{upward directed}, if every pair of elements of~$P$ is contained in~$P\dnw x$ for some~$x\in P$. We say that~$U$ is \emph{cofinal} in~$P$, if $P\dnw U=P$. An \emph{ideal} of~$P$ is a nonempty, upward directed, lower subset of~$P$. We set
 \[
 P^{[2]}:=\setm{(x,y)\in P\times P}{x\leq y}\,.
 \]
For partially ordered sets~$P$ and~$Q$, a map $f\colon P\to Q$ is \emph{isotone} (resp., \emph{antitone}), if $x\leq y$ implies that $f(x)\leq f(y)$ (resp., $f(y)\leq f(x)$), for all $x,y\in P$.

We refer to Birkhoff~\cite{Birk94} or Gr\"atzer~\cite{GLT2} for basic notions of lattice theory. We recall here a sample of needed notation, terminology, and results. A family $\famm{a_i}{i\in I}$ of elements in a lattice~$L$ with zero is \emph{independent}, if the equality
 \[
 \bigvee\famm{a_i}{i\in X}\wedge\bigvee\famm{a_i}{i\in Y}=
 \bigvee\famm{a_i}{i\in X\cap Y}
 \]
holds for all finite subsets~$X$ and~$Y$ of~$I$. In case~$L$ is modular and~$I=\set{0,\dots,n-1}$ for a non-negative integer~$n$, this amounts to verifying that the equality $a_k\wedge\bigvee_{i<k}a_i=0$ holds for each $k<n$. We denote by~$\oplus$ the operation of finite independent sum in~$L$; hence $a=\bigoplus\famm{a_i}{i\in I}$ means that~$I$ is finite, $\famm{a_i}{i\in I}$ is independent, and $a=\bigvee_{i<n}a_i$. If~$L$ is modular, then~$\oplus$ is both commutative and associative in the strongest possible sense for a partial operation, see \cite[Section~II.1]{Maed58}.

A lattice~$L$ with zero is \emph{sectionally complemented}, if for all $a\leq b$ in~$L$ there exists~$x\in L$ such that $b=a\oplus x$. For elements $a,x,b\in L$, let $a\sim_xb$ hold, if $a\oplus x=b\oplus x$. We say that~$a$ is \emph{perspective} to~$b$, in notation~$a\sim b$, if there exists~$x\in L$ such that $a\sim_xb$.
We say that~$L$ is \emph{complemented}, if it has a unit and every element~$a\in L$ has a \emph{complement}, that is, an element~$x\in L$ such that $1=a\oplus x$. A bounded modular lattice is complemented if and only if it is sectionally complemented.

An ideal~$I$ of a lattice~$L$ is \emph{neutral}, if $\set{I,X,Y}$ generates a distributive sublattice of~$\Id L$ for all ideals~$X$ and~$Y$ of~$L$. In case~$L$ is sectionally complemented modular, this is equivalent to the statement that every element of~$L$ perspective to some element of~$I$ belongs to~$I$. In that case, the assignment that to a congruence~$\theta$ associates the $\theta$-block of~$0$ is an isomorphism from the congruence lattice of~$L$ onto the lattice of all neutral ideals of~$L$.

An independent finite sequence $\famm{a_i}{i<n}$ in a lattice~$L$ with zero is \emph{homogeneous}, if the elements~$a_i$ are pairwise perspective. An element~$x\in L$ is \emph{large}, if the neutral ideal generated by~$x$ is~$L$.

A pair $(\famm{a_i}{0\leq i<n},\famm{c_i}{1\leq i<n})$, with $\famm{a_i}{0\leq i<n}$ independent, is a
\begin{itemize}
\item\emph{$n$-frame}, if $a_0\sim_{c_i}a_i$ for each~$i$ with $1\leq i<n$;

\item\emph{large $n$-frame}, if it is an $n$-frame and~$a_0$ is large.
\end{itemize}

The assignment $R\mapsto\LL(R)$ extends canonically to a \emph{functor} from the category of all regular rings with ring homomorphisms to the category of \scml s with $0$-lattice homomorphisms (cf. Micol \cite{Micol} for details). This functor preserves direct limits.

Denote by $\Idemp R$ the set of all idempotent elements in a ring~$R$.
For idempotents $a$ and $b$ in a ring~$R$, let $a\utr b$ hold, if $a=ab=ba$; equivalently, $a\in bRb$.

We shall need the following folklore lemma.

\begin{lemma}\label{L:DecompId}
Let $A$ and $B$ be right ideals in a ring~$R$ and let~$e$ be an idempotent element of~$R$. If $eR=A\oplus B$, then there exists a unique pair $(a,b)\in A\times B$ such that $e=a+b$. Furthermore, both~$a$ and~$b$ are idempotent, $ab=ba=0$, $A=aR$, and $B=bR$.
\end{lemma}

\section{\Ban\ functions on lattices and rings}\label{S:BanLattRing}

\begin{definition}\label{D:BanLatt}
Let $X$ be a subset in a bounded lattice~$L$. A \emph{partial \Ban\ function on~$X$ in~$L$} is an antitone map $f\colon X\to L$ such that $x\oplus f(x)=1$ for each $x\in X$. In case~$X=L$, we say that~$f$ is a \emph{\Ban\ function on~$L$}.
\end{definition}

Trivially, every bounded lattice with a \Ban\ function is complemented. The following example shows that the converse does not hold as a rule.

\begin{example}\label{Ex:NonBanCL}
The finite lattice~$F$ diagrammed on Figure~\ref{Fig:NonBan} is complemented. However, $F$ does not have any \Ban\ function, because~$a'$ is the unique complement of~$a$, $b'$ is the unique complement of~$b$, $a\leq b$, while $b'\nleq a'$.
\begin{figure}[htb]
\includegraphics{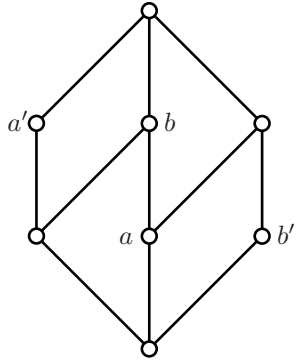}
\caption{A finite complemented lattice without a \Ban\ function}
\label{Fig:NonBan}
\end{figure}

\end{example}

Although most lattices involved in the present paper will be modular, it is noteworthy to observe that \Ban\ functions may also be of interest in the `orthogonal' case of \msd\ lattices. By definition, a lattice~$L$ is \emph{\msd}, if $x\wedge y=x\wedge z$ implies that $x\wedge y=x\wedge(y\vee z)$, for all $x,y,z\in L$. The following result has been pointed to the author by Luigi Santocanale.

\begin{proposition}\label{P:Luigi}
Let $L$ be finite lattice. Consider the following conditions:
\begin{enumerate}
\item the set of all atoms of~$L$ joins to the largest element of~$L$;

\item $L$ has a \Ban\ function;

\item $L$ is complemented.
\end{enumerate}
Then \textup{(ii)} implies \textup{(iii)} implies \textup{(i)}. Furthermore, if~$L$ is \msd, then~\textup{(i)}, \textup{(ii)}, and \textup{(iii)} are equivalent.
\end{proposition}

\begin{proof}
Denote by~$\At L$ the set of all atoms of~$L$.

(i)$\Rightarrow$(ii) \emph{in case~$L$ is \msd}. Set
 \[
 f(x):=\bigvee\famm{p\in\At L}{p\wedge x=0}\,,
 \]
for each $x\in L$. For $x\in L$ and $p\in\At L$, if $p\nleq x\vee f(x)$, then $p\nleq x$, thus, as~$p$ is an atom, $p\wedge x=0$, thus, by the definition of~$f$, $p\leq f(x)$, a contradiction. Thus $p\leq x\vee f(x)$ for each $p\in\At L$, and thus, by assumption, $x\vee f(x)=1$. Furthermore, it follows from the \msdy\ of~$L$ that $x\wedge f(x)=0$, for each $x\in L$. As~$f$ is obviously antitone, $f$ is a \Ban\ function on~$L$.

(ii)$\Rightarrow$(iii) is trivial.

(iii)$\Rightarrow$(i). Set $a:=\bigvee\At L$. As~$L$ is complemented, there exists $b\in L$ such that $a\oplus b=1$. If~$b$ is nonzero, then there exists an atom~$p$ below~$b$, thus $p\nleq a$, a contradiction. Hence $b=0$, and so $a=1$.
\end{proof}

The conditions (i)--(iii) of Proposition~\ref{P:Luigi} are not uncommon. They are, for example, satisfied for the permutohedron on a given finite number of letters. It follows that they are also satisfied for the associahedron (Tamari lattice), which is a quotient of the permutohedron.

We shall now introduce a ring-theoretical analogue of the definition of a \Ban\ function.

\begin{definition}\label{D:BanRing}
Let $X$ be a subset in a ring~$R$. A \emph{partial \Ban\ function on~$X$ in~$R$} is a mapping $f\colon X\to\Idemp R$ such that
\begin{enumerate}
\item $xR=f(x)R$ for each $x\in X$.

\item $xR\subseteq yR$ implies that $f(x)\utr f(y)$, for all $x,y\in X$.
\end{enumerate}
In case $X=R$ we say that~$f$ is a \emph{\Ban\ function on~$R$}.
\end{definition}

In the context of Definition~\ref{D:BanRing}, we put
 \begin{equation}\label{Eq:DefnLR(X)}
 \LL_R(X):=\setm{xR}{x\in X}\,.
 \end{equation}
\Ban\ functions in rings and in lattices are related by the following result.

\begin{lemma}\label{L:BanLattRing}
Let~$R$ be a unital von~Neumann regular ring and let $X\subseteq R$. Then the following are equivalent:
\begin{enumerate}
\item There exists a partial \Ban\ function on~$\LL_R(X)$ in~$\LL(R)$.

\item There exists a partial \Ban\ function on~$X$ in~$R$.
\end{enumerate}
\end{lemma}

\begin{proof}
(i)$\Rightarrow$(ii). Let $\varphi\colon\LL_R(X)\to\LL(R)$ be a partial \Ban\ function. For each $x\in X$, as $R=xR\oplus\varphi(xR)$ it follows from Lemma~\ref{L:DecompId} that the unique element $f(x)\in xR$ such that $1-f(x)\in\varphi(xR)$ is idempotent and satisfies both relations $xR=f(x)R$ and $\varphi(xR)=(1-f(x))R$. Let $x,y\in X$ such that $xR\subseteq yR$. {}From $f(x)R=xR\subseteq yR=f(y)R$ and the idempotence of~$f(y)$ it follows that $f(x)=f(y)f(x)$. {}From
 \[
 (1-f(y))R=\varphi(yR)\subseteq\varphi(xR)=(1-f(x))R
 \]
together with the idempotence of~$f(x)$ we get $f(x)(1-f(y))=0$, and thus $f(x)=\nobreak f(x)f(y)$. Therefore, $f(x)\utr f(y)$.

(ii)$\Rightarrow$(i). Let $f\colon X\to\Idemp R$ be a partial \Ban\ function. As
 \[
 xR\subseteq yR\Rightarrow f(x)\utr f(y)\Rightarrow 1-f(y)\utr 1-f(x)
 \Rightarrow(1-f(y))R\subseteq(1-f(x))R\,,
 \]
there exists a unique map $\varphi\colon\LL_R(X)\to\LL(R)$ such that
 \[
 \varphi(xR)=(1-f(x))R\,,\quad\text{for each }x\in X\,,
 \]
and $\varphi$ is antitone. Furthermore, for each $x\in X$, from the idempotence of~$f(x)$ it follows that $R=f(x)R\oplus(1-f(x))R$, that is, $R=xR\oplus\varphi(xR)$. Therefore, $\varphi$ is a partial \Ban\ function on~$\LL_R(X)$ in~$\LL(R)$.
\end{proof}

\section{\Ban\ functions on countable \cml s}\label{S:CML}

A large part of the present section will be devoted to proving the following result.

\begin{theorem}\label{T:BaCML}
Every countable \cml\ has a \Ban\ function with Boolean range.
\end{theorem}

Let~$L$ be a \cml. We denote by~$\cB$ the set of all finite sequences $\bdu =\famm{u_i}{i<n}$, where $n=:|\bdu|<\omega$, of elements of~$L$ such that $1=\bigoplus\famm{u_i}{i<n}$. We set $Z(\bdu):=\setm{k<|\bdu|}{u_k=0}$, and, further, $u_{<k}:=\bigvee\famm{u_i}{i<k}$ for each $k\leq|\bdu|$ (with $u_{<0}:=0$). Furthermore, for each $x\in L$ we set
 \begin{align*}
 F_\bdu(x)&:=\setm{k<|\bdu|}{u_k\nleq x\vee u_{<k}}\,,\\
 G_\bdu(x)&:=\setm{k<|\bdu|}{u_k\wedge(x\vee u_{<k})=0}\,,\\ 
 f_\bdu(x)&:=\bigvee\famm{u_k}{k\in F_\bdu(x)}\,,\\
 g_\bdu(x)&:=\bigvee\famm{u_k}{k\in G_\bdu(x)}\,.
 \end{align*}

\begin{lemma}\label{L:xvfux=1}
The following statements hold, for each $\bdu\in\cB$ and each $x\in L$:
\begin{enumerate}
\item $x\vee f_\bdu(x)=1$;

\item $x\wedge g_\bdu(x)=0$;

\item $g_\bdu(x)\leq f_\bdu(x)$.
\end{enumerate}
\end{lemma}

\begin{proof}
(i). As $\bigvee\famm{u_k}{k<|\bdu|}=1$, it suffices to prove that $u_k\leq x\vee f_\bdu(x)$ for each $k<\nobreak|\bdu|$. We argue by induction on~$k$; the induction hypothesis is that\linebreak $u_{<k}\leq x\vee f_\bdu(x)$. If $u_k\leq x\vee u_{<k}$ then, by the induction hypothesis, $u_k\leq x\vee f_\bdu(x)$ as well, while if $u_k\nleq x\vee u_{<k}$, that is, $k\in F_\bdu(x)$, then $u_k\leq f_\bdu(x)\leq x\vee f_\bdu(x)$.

(ii). For each $k\in G_\bdu(x)$, from $u_k\wedge(x\vee u_{<k})=0$ it follows \emph{a fortiori} that $u_k\wedge\bigl(x\vee\bigvee\famm{u_i}{i<k\text{ and }i\in G_\bdu(x)}\bigr)=0$. Therefore, writing $G_\bdu(x)=\setm{k_s}{s<r}$ with $k_0<\cdots<k_{r-1}$, we obtain, by using the modularity of~$L$, that the finite sequence $(x,u_{k_0},\dots,u_{k_{r-1}})$ is independent in~$L$. In particular,
 \begin{equation*}
 x\wedge g_\bdu(x)=x\wedge\bigvee\famm{u_{k_s}}{s<r}=0\,.\tag*{\qed}
 \end{equation*}
(iii) follows immediately from the containment $G_\bdu(x)\subseteq F_\bdu(x)\cup Z(\bdu)$.
\end{proof}

\begin{lemma}\label{L:fuxpartBan}
Let $\bdu\in\cB$ and let $x,y\in L$. If $x\leq y$, then $f_\bdu(y)\leq f_\bdu(x)$.
\end{lemma}

\begin{proof}
{}From the inequality $x\leq y$ it follows that $F_\bdu(y)\subseteq F_\bdu(x)$. The conclusion follows immediately from the definition of~$f_\bdu$.
\end{proof}

For $\bdu,\bdv\in\cB$ and $\varphi\colon\set{0,\dots,|\bdv|-1}\onto\set{0,\dots,|\bdu|-1}$ isotone and surjective, let $\varphi\colon\bdv\onto\bdu$ hold, if
 \begin{equation}\label{Eq:ContrMorcB}
 u_k=\bigvee\famm{v_l}{l\in\varphi^{-1}\set{k}}\quad\text{for each }k<|\bdu|
 \end{equation}
(observe that the join in~\eqref{Eq:ContrMorcB} is necessarily independent). We say that~$\bdv$ \emph{refines}~$\bdu$, if there exists~$\varphi$ such that $\varphi\colon \bdv\onto \bdu$. Then we denote by $\varphi_-(k)$ (resp., $\varphi_+(k)$) the least (resp., largest) element of $\varphi^{-1}\set{k}$, for each $k<|\bdu|$. As~$\varphi$ is isotone and surjective, $\varphi_-(k)\leq\varphi_+(k)$ and $\varphi^{-1}\set{k}=[\varphi_-(k),\varphi_+(k)]$.

Say that an element $\bdu\in\cB$ \emph{decides} an element $x\in L$, if $F_\bdu(x)\subseteq G_\bdu(x)$. By Lemma~\ref{L:xvfux=1}(iii), it follows that $f_\bdu(x)=g_\bdu(x)$.

\begin{lemma}\label{L:Ineqfugu}
Let $\bdu,\bdv\in\cB$, let $\varphi\colon\bdv\onto\bdu$, and let $x\in L$. Then the following statements hold:
\begin{enumerate}
\item $v_l\leq u_{\varphi(l)}$ and $u_{<\varphi(l)}\leq v_{<l}$, for each $l<|\bdv|$.

\item $\varphi F_\bdv(x)\subseteq F_\bdu(x)$;

\item $\varphi^{-1}G_\bdu(x)\subseteq G_\bdv(x)$;

\item $f_\bdv(x)\leq f_\bdu(x)$;

\item $g_\bdu(x)\leq g_\bdv(x)$;

\item if $\bdv$ refines~$\bdu$ and~$\bdu$ decides~$x$, then~$\bdv$ decides~$x$ and $f_\bdu(x)=f_\bdv(x)$.
\end{enumerate}
\end{lemma}

\begin{proof}
(i) follows easily from~\eqref{Eq:ContrMorcB}.

(ii). Let $l\in F_\bdv(x)$ and set $k:=\varphi(l)$. {}From $v_l\nleq x\vee v_{<l}$ together with~(i) it follows that $u_k\nleq x\vee u_{<k}$, that is, $k\in F_\bdu(x)$.

(iii). Let $l\in\varphi^{-1}G_\bdu(x)$, so $k:=\varphi(l)$ belongs to $G_\bdu(x)$, that is, $u_k\wedge(x\vee u_{<k})=0$. As~$L$ is modular and by~\eqref{Eq:ContrMorcB}, this means that the finite sequence
 \[
 (x\vee u_{<k},v_{\varphi_-(k)},\dots,v_{\varphi_+(k)})
 \]
is independent, thus, as $\varphi_-(k)\leq l\leq\varphi_+(k)$,
 \[
 v_l\wedge\bigl(x\vee u_{<k}\vee\bigvee\famm{v_i}{\varphi_-(k)\leq i<l}\bigr)=0\,,
 \]
that is, by~\eqref{Eq:ContrMorcB}, $v_l\wedge(x\vee v_{<l})=0$, which means that $l\in G_\bdv(x)$.

(iv). For each $l\in F_\bdv(x)$, it follows from~(i) that $v_l\leq u_{\varphi(l)}$ and from~(ii) that $\varphi(l)\in F_\bdu(x)$, thus $v_l\leq u_{\varphi(l)}\leq f_\bdu(x)$. As this holds for each $l\in F_\bdv(x)$, we obtain that $f_\bdv(x)\leq f_\bdu(x)$.

(v). Let $k\in G_\bdu(x)$. It follows from~(iii) that $\varphi^{-1}\set{k}\subseteq G_\bdv(x)$, thus, by~\eqref{Eq:ContrMorcB}, $u_k\leq\bigvee\famm{v_l}{l\in G_\bdv(x)}=g_\bdv(x)$. This holds for each $k\in G_\bdu(x)$, thus $g_\bdu(x)\leq g_\bdv(x)$.

(vi). As $F_\bdu(x)\subseteq G_\bdu(x)$, we obtain, by using~(ii) and~(iii),
 \[
 F_\bdv(x)\subseteq\varphi^{-1}\varphi F_\bdv(x)\subseteq\varphi^{-1} F_\bdu(x)
 \subseteq\varphi^{-1}G_\bdu(x)\subseteq G_\bdv(x)\,,
 \]
so~$\bdv$ decides~$x$. As both~$\bdu$ and~$\bdv$ decide~$x$, we obtain that $f_\bdu(x)=g_\bdu(x)$ and $f_\bdv(x)=g_\bdv(x)$, so the conclusion follows from~(iv) and~(v).
\end{proof}

\begin{lemma}\label{L:forcvdecx}
For each $\bdu\in\cB$ and each $x\in L$, there exists $\bdv\in\cB$ such that~$\bdv$ refines~$\bdu$ and~$\bdv$ decides~$x$.
\end{lemma}

\begin{proof}
Set $n:=|\bdu|$. For each $k<n$, we set $v_{2k}:=u_k\wedge(x\vee u_{<k})$ and we pick~$v_{2k+1}$ such that $u_k=v_{2k}\oplus v_{2k+1}$. It is obvious that the finite sequence $v:=\famm{v_l}{l<2n}$ belongs to~$\cB$ and refines~$\bdu$.

It remains to verify that~$\bdv$ decides~$x$. So let $l<2n$. If $l=2k$ for some $k<n$, then $v_l=v_{2k}\leq x\vee u_{<k}=x\vee v_{<l}$. Suppose that $l=2k+1$ for some $k<n$. As $u_i=v_{2i}\vee v_{2i+1}$ for each $i<k$ while $v_{2k}\leq x\vee u_{<k}$, we get
 \[
 x\vee v_{<l}=x\vee v_0\vee v_1\vee\cdots\vee v_{2k}
 =x\vee u_{<k}\vee v_{2k}=x\vee u_{<k}\,,
 \]
so $v_l\wedge(x\vee v_{<l})=v_{2k+1}\wedge(x\vee u_{<k})
=v_{2k+1}\wedge u_k\wedge(x\vee u_{<k})=v_{2k+1}\wedge v_{2k}=0$.
\end{proof}

\begin{proof}[Proof of Theorem~\textup{\ref{T:BaCML}}]
As~$L$ is countable, we can write $L=\setm{a_n}{n<\omega}$ and denote by~$\nu(x)$ the least non-negative integer~$n$ such that $x=a_n$, for each $x\in L$. It follows from Lemmas~\ref{L:Ineqfugu}(vi) and~\ref{L:forcvdecx} that there exists a sequence $\famm{\bdu_n}{n<\omega}$ of elements of~$\cB$ such that~$\bdu_n$ decides all elements $a_0$, \dots, $a_n$ and~$\bdu_{n+1}$ refines~$\bdu_n$, for each $n<\omega$. We set $f(x):=f_{\bdu_{\nu(x)}}(x)$, for each $x\in L$. Observe that, by Lemma~\ref{L:Ineqfugu}(vi), $f(x)=f_{\bdu_n}(x)$ for each integer $n\geq\nu(x)$. {}From Lemma~\ref{L:xvfux=1} it follows that $1=x\oplus f(x)$. Finally, from Lemma~\ref{L:fuxpartBan} it follows that the map~$f$ is antitone, so it is a \Ban\ function on~$L$.

Furthermore, (the underlying set of) each~$\bdu_n$ is independent with join~$1$, thus it generates a Boolean sublattice~$B_n$ of~$L$ with the same bounds as~$L$. As~$\bdu_{n+1}$ refines~$\bdu_n$, $B_{n+1}$ contains~$B_n$. As the range of each~$f_{\bdu_n}$ is contained in~$B_n$, the range of~$f$ is contained in the Boolean sublattice $B:=\bigcup\famm{B_n}{n<\omega}$ of~$L$. For each $x\in B$, $f(x)$ is a complement of~$x$ in~$B$, thus it is the unique complement of~$x$ in~$B$---denote it by~$\neg x$. As $B=\setm{\neg x}{x\in B}$, it follows that the range of~$f$ is exactly~$B$.
\end{proof}

For von~Neumann regular rings we get the following corollary.

\begin{corollary}\label{C:BanRegRings}
Every countable von~Neumann regular ring has a \Ban\ function.
\end{corollary}

We emphasize that we do not require the ring be unital in Corollary~\ref{C:BanRegRings}.

\begin{proof}
Let~$R$ be a countable von~Neumann regular ring.
By Fuchs and Hal\-pe\-rin~\cite{FuHa64}, $R$ embeds as a two-sided ideal into some \emph{unital} von~Neumann regular ring~$S$. Starting with~$R\cup\set{1}$ and closing under the ring operations and a given operation of quasi-inversion on~$S$, we obtain a countable von~Neumann regular subring of~$S$ containing~$R\cup\set{1}$; hence we may assume that~$S$ is countable. It follows from Theorem~\ref{T:BaCML} that~$\LL(S)$ has a \Ban\ function. By Lemma~\ref{L:BanLattRing}, it follows that~$S$ has a \Ban\ function, say~$g$. For each~$x\in R$, as $xS=g(x)S$ and~$R$ is a right ideal of~$S$, $g(x)$ belongs to~$R$. Furthermore, there exists $y\in S$ such that $g(x)=xy$, thus, as~$g(x)$ is idempotent, $g(x)=xyxy$. As~$R$ is a two-sided ideal of~$S$, $yxy$ belongs to~$R$, and thus~$g(x)$ belongs to~$xR$. As~$x=g(x)x$, it follows that $xR=g(x)R$. It follows that the restriction of~$g$ from~$R$ to $\Idemp R$ is a \Ban\ function on~$R$.
\end{proof}

Say that a \Ban\ function on a lattice~$L$ is \emph{Boolean}, if its range is a Boolean sublattice of~$L$. 
In case~$L$ is the subspace lattice of a vector space~$V$, the range~$B$ of a Boolean \Ban\ function on~$L$ may be chosen as the set of all spans of all subsets of a given basis of~$V$. In particular, $B$ is far from being unique.

However, we shall now prove that if~$L$ is a countable \cml, then~$B$ is unique \emph{up to isomorphism}. For a Boolean algebra~$B$ and a \cm~$M$, a \emph{V-measure} (cf. Dobbertin~\cite{Dobb83}) from~$B$ to~$M$ is a map $\mu\colon B\to M$ such that $\mu(x)=0$ if and only if $x=0$, $\mu(x\oplus y)=\mu(x)+\mu(y)$ for all disjoint $x,y\in B$, and if $\mu(z)=\alpha+\beta$, then there are $x,y\in B$ such that $z=x\oplus y$, $\mu(x)=\alpha$, and $\mu(y)=\beta$.

Denote by~$\Delta$ the canonical map from~$L$ to its \emph{dimension monoid}~$\Dim L$, see page~259 and Chapter~9 in Wehrung~\cite{WDim}. 

\begin{proposition}\label{P:DeltaBVm}
Let $f$ be a \Ban\ function with Boolean range~$B$ on a \cml\ $L$. Then the restriction of~$\Delta$ from~$B$ to~$\Dim L$ is a V-measure on~$B$.
\end{proposition}

\begin{proof}
It is obvious that $\Delta(x)=0$ if and only if $x=0$, for each $x\in L$, and that $\Delta(x\vee y)=\Delta(x)+\Delta(y)$ whenever~$x$ and~$y$ are disjoint elements in~$B$ (for they are also disjoint in~$L$). Now let $c\in B$ and let $\alpha,\beta\in\Dim L$ such that $\Delta(c)=\alpha+\beta$. It follows from \cite[Corollary~9.6]{WDim} that there are $x,y\in L$ such that $c=x\oplus y$, $\Delta(x)=\alpha$, and $\Delta(y)=\beta$.

Put $b:=c\wedge f(x)$. As both~$c$ and~$f(x)$ belong to~$B$, the element~$b$ also belongs to~$B$. Furthermore, $x\wedge b=x\wedge f(x)=0$, and
 \begin{align*}
 c&=c\wedge(x\vee f(x))\\
 &=x\vee(c\wedge f(x))
 &&(\text{because }x\leq c\text{ and }L\text{ is modular})\\
 &=x\vee b\,, 
 \end{align*}
so $c=x\oplus y=x\oplus b$ and so~$y$ and~$b$ are perspective. In particular, $\Delta(b)=\Delta(y)=\beta$.

Likewise, there exists $a\in B$ such that $c=x\oplus b=a\oplus b$, so $\Delta(a)=\Delta(x)=\alpha$.
\end{proof}

For Boolean algebras~$A$ and~$B$, a subset~$\rho$ of~$A\times B$ is an \emph{additive V-relation}, if $1_A\mathbin{\rho}1_B$, $x\mathbin{\rho}0_B$ if and only if $x=0_A$, $x\mathbin{\rho}y_0\oplus y_1$ if and only if there exists a decomposition $x=x_0\oplus x_1$ with $x_0\mathbin{\rho}y_0$ and $x_1\mathbin{\rho}y_1$, and symmetrically with~$A$ and~$B$ interchanged.
Vaught's isomorphism Theorem (cf. \cite[Theorem~1.1.3]{Pier}) implies that any additive V-relation between countable Boolean algebras~$A$ and~$B$ contains the graph of some isomorphism from~$A$ onto~$B$.

In particular, if~$A$ and~$B$ are Boolean algebras, then, for any V-measures $\lambda\colon A\to\nobreak M$ and $\mu\colon B\to M$ such that $\lambda(1_A)=\mu(1_B)$, the binary relation
 \[
 R:=\setm{(x,y)\in A\times B}{\lambda(x)=\mu(y)}
 \]
is an additive V-relation between~$A$ and~$B$. Therefore, if both~$A$ and~$B$ are countable, then, by Vaught's Theorem, there exists an isomorphism $\varphi\colon A\to B$ such that $\lambda=\mu\circ\varphi$.

By the above paragraph, we obtain

\begin{corollary}\label{C:DeltaBVm}
Let~$L$ be a countable \cml. Then for a Boolean Banaschewski function on~$L$ with range~$B$, the pair $(B,\Delta\res_B)$ is unique up to isomorphism. In particular, $B$ is unique up to isomorphism.
\end{corollary}

\section{\Ban\ measures and \Ban\ traces}\label{S:BanTail}

\begin{definition}\label{D:BanTail}
A \emph{\Ban\ trace} on a lattice~$L$ with zero is a family\linebreak $\famm{a_i^j}{i\leq j\text{ in }\Lambda}$ of elements in~$L$, where~$\Lambda$ is an upward directed partially ordered set with zero, such that
\begin{enumerate}
\item $a_i^k=a_i^j\oplus a_j^k$ for all $i\leq j\leq k$ in~$\Lambda$;

\item $\setm{a_0^i}{i\in\Lambda}$ is cofinal in~$L$.
\end{enumerate}
We say that the \Ban\ trace above is \emph{normal}, if $i\leq j$ and $a_0^i=a_0^j$ implies that $i=j$, for all $i,j\in\Lambda$.
\end{definition}

It is trivial that every bounded lattice has a normal \Ban\ trace (if $L=\set{0}$ take $\Lambda=\set{0}$ and $a_0^0=0$; if~$L$ is bounded nontrivial take $\Lambda=\set{0,1}$ and $a_0^0=a_1^1=0$ while $a_0^1=1$), so this notion is interesting only for unbounded lattices.

It is obvious that every \scml\ embeds into a reduced product of its principal ideals, thus into a \cml. Our first application of \Ban\ traces, namely Theorem~\ref{T:BanTailId}, deals with the question whether such an embedding can be taken with ideal range. We will use the following well-known lemma.

\begin{lemma}[Folklore]\label{L:xyzindep}
Let $x$, $y$, $z$ be elements in a modular lattice~$L$. If\linebreak $(x\vee y)\wedge z\leq y$, then $x\wedge(y\vee z)=x\wedge y$ and $(x\vee z)\wedge(y\vee z)=(x\wedge y)\vee z$.
\end{lemma}

\begin{note}
It is not hard to verify that the conclusion of Lemma~\ref{L:xyzindep} can be strengthened by stating that the sublattice of~$L$ generated by $\set{x,y,z}$ is \emph{distributive}.
\end{note}

\begin{proof}
We start by computing, using the modularity of~$L$ and the assumption,
 \[
 (x\vee y)\wedge(y\vee z)=y\vee\bigl((x\vee y)\wedge z)\bigr)=y\,.
 \]
It follows that
 \[
 x\wedge(y\vee z)=x\wedge(x\vee y)\wedge(y\vee z)=x\wedge y\,.
 \]
It follows, by using again the modularity of~$L$, that
 \begin{equation*}
 (x\vee z)\wedge(y\vee z)=\bigl(x\wedge(y\vee z)\bigr)\vee z=(x\wedge y)\vee z\,.
 \tag*{\qed}
 \end{equation*}
\renewcommand{\qed}{}
\end{proof}

\begin{theorem}\label{T:BanTailId}
Every \scml\ with a \Ban\ trace embeds, as a neutral ideal and within the same quasivariety, into some \cml.
\end{theorem}

\begin{proof}
Let $\famm{a_i^j}{i\leq j\text{ in }\Lambda}$ be a \Ban\ trace in a \scml~$L$. The conclusion of the theorem for~$L$ is trivial in case~$L$ has a unit, so suppose that~$L$ has no unit.

We denote by~$\cF$ the filter on~$\Lambda$ generated by all principal upper subsets $\Lambda\upw i$, for $i\in\Lambda$, and we denote by~$\ol{L}$ the reduced product of the family $\famm{L\dnw a_0^i}{i\in\Lambda}$ modulo~$\cF$. For any $i_0\in\Lambda$ and any family $\famm{x_i}{i\in\Lambda\upw i_0}$ in $\prod_{i\in\Lambda\upw i_0}(L\dnw a_0^i)$, we shall denote by $[x_i\mid i\to\infty]$ the equivalence class modulo~$\cF$ of the family $\famm{y_i}{i\in\Lambda}$ defined by
 \[
 y_i:=\begin{cases}
 x_i\,,&\text{if }i\geq i_0\,,\\
 0\,,&\text{otherwise},
 \end{cases}
 \quad\text{for every }i\in\Lambda\,.
 \]
In particular, for each $x\in L$, the subset $\setm{j\in\Lambda}{x\leq a_0^j}$ contains a principal filter of~$\Lambda$, thus we can define a map $\eps\colon L\to\ol{L}$ by the rule
 \[
 \eps(x):=[x\mid j\to\infty]\,,\quad\text{for each }x\in L\,.
 \]
Furthermore, for each $i\in\Lambda$, define a map $\eps_i\colon L\dnw a_0^i\to\ol{L}$ by the rule
 \[
 \eps_i(x):=[x\vee a_i^j\mid j\to\infty]\,,\quad\text{for each }x\in L\dnw a_0^i\,.
 \]
Consider the following subset of~$\ol{L}$.
 \begin{equation}\label{Eq:DefntildeL}
 \tilde{L}:=\im\eps\cup\bigcup\famm{\im\eps_i}{i\in\Lambda}\,.
 \end{equation}
The following claim shows that the union on the right hand side of~\eqref {Eq:DefntildeL} is directed.

\setcounter{claim}{0}
\begin{claim}\label{Cl:epsileqepsj}
$i\leq j$ implies that $\im\eps_i\subseteq\im\eps_j$, for all $i,j\in\Lambda$.
\end{claim}

\begin{cproof}
For all $x\in L\dnw a_0^i$,
 \begin{align*}
 \eps_i(x)&=[x\vee a_i^k\mid k\to\infty]\\
 &=[x\vee a_i^j\vee a_j^k\mid k\to\infty]\\
 &=\eps_j(x\vee a_i^j)\,.\tag*{\qedc} 
 \end{align*}
\renewcommand{\qedc}{}
\end{cproof}
Now it is obvious that~$\eps$ is a $0$-lattice embedding from~$L$ into~$\ol{L}$, while~$\eps_i$ is a \jh, for each $i\in\Lambda$. Furthermore, $\eps(x)\vee\eps_i(y)=\eps_i(x\vee y)$, for all $i\in\Lambda$ and all $x,y\in L\dnw a_0^i$. In particular, by Claim~\ref{Cl:epsileqepsj}, the subset~$\tilde{L}$ defined in~\eqref{Eq:DefntildeL} is a \jz-subsemilattice of~$\ol{L}$.

\begin{claim}\label{Cl:epsialmhom}
Let $i\in\Lambda$ and let $x,y\in L\dnw a_0^i$.
Then both equalities $\eps(x)\wedge\eps_i(y)=\eps(x\wedge y)$ and $\eps_i(x)\wedge\eps_i(y)=\eps_i(x\wedge y)$ hold. In particular, $\eps_i$ is a lattice homomorphism from~$L\dnw a_0^i$ to~$\ol{L}$.
\end{claim}

\begin{cproof}
Let $j\in\Lambda\upw i$. {}From $x\vee y\leq a_0^i$ and $a_0^i\wedge a_i^j=0$ it follows that $(x\vee y)\wedge a_i^j=0$. By Lemma~\ref{L:xyzindep}, we obtain the following equations:
 \[
 x\wedge(y\vee a_i^j)=x\wedge y\quad\text{and}\quad
 (x\vee a_i^j)\wedge(y\vee a_i^j)=(x\wedge y)\vee a_i^j\,.
 \]
Therefore, by evaluating the equivalence class modulo~$\cF$ of both sides of each of the equalities above as $j\to\infty$, we obtain the desired conclusion.
\end{cproof}

In particular, from Claims~\ref{Cl:epsileqepsj} and~\ref{Cl:epsialmhom} it follows that~$\tilde{L}$ is a meet-subsemilattice of~$\ol{L}$. Therefore, \emph{$\tilde{L}$ is a $0$-sublattice of~$\ol{L}$}. As~$\ol{L}$ is a reduced product of sublattices of~$L$, it belongs to the same quasivariety as~$L$; hence so does~$\tilde{L}$.

Furthermore, for all $x,y\in L$ and all $i\in\Lambda$ such that $x\vee y\leq a_0^i$, if $\eps_i(y)\leq\eps(x)$, then, by Claim~\ref{Cl:epsialmhom},
 \[
 \eps_i(y)=\eps_i(y)\wedge\eps(x)=\eps(x\wedge y)\,,
 \]
thus~$\eps_i(y)$ belongs to~$\im\eps$. Therefore, \emph{$\im\eps$ is an ideal of~$\tilde{L}$}.

Now we verify that $\tilde{L}$ is a \cml. It has a unit, namely $1_{\tilde{L}}=\eps_0(0)=[a_0^i\mid i\to\infty]$. Let $x\in L$ and let $i\in\Lambda$ such that $x\leq a_0^i$. As~$L$ is sectionally complemented, there exists $y\in L\dnw a_0^i$ such that $x\oplus y=a_0^i$. Hence
 \[
 \eps(x)\vee\eps_i(y)=\eps_i(x\vee y)=\eps_i(a_0^i)=
 [a_0^i\vee a_i^j\mid j\to\infty]=[a_0^j\mid j\to\infty]=1_{\tilde{L}}\,,
 \]
while, by Claim~\ref{Cl:epsialmhom},
 \[
 \eps(x)\wedge\eps_i(y)=\eps(x\wedge y)=\eps(0)=0\,.
 \]
Therefore, $1_{\tilde{L}}=\eps(x)\oplus\eps_i(y)$. By symmetry between~$x$ and~$y$, we also obtain $1_{\tilde{L}}=\eps_i(x)\oplus\eps(y)$. Therefore, $\tilde{L}$ is complemented.

It remains to prove that~$\im\eps$ is a neutral ideal of~$\tilde{L}$. By \cite[Theorem~III.20]{Birk94}, it suffices to prove that~$\im\eps$ contains any element of~$\tilde{L}$ perspective to some element of~$\im\eps$. By using Claim~\ref{Cl:epsileqepsj}, it suffices to prove that for any $i\in\Lambda$ and any $x,y,z\in L\dnw a_0^i$, none of the relations $\eps_i(x)\sim_{\eps(z)}\eps(y)$ and $\eps_i(x)\sim_{\eps_i(z)}\eps(y)$ can occur.

If $\eps_i(x)\sim_{\eps(z)}\eps(y)$, then $\eps_i(x\vee z)=\eps_i(x)\vee\eps(z)=\eps(y)\vee\eps(z)=\eps(y\vee z)$, thus there exists $j\in\Lambda\upw i$ such that
 \[
 x\vee z\vee a_i^k=y\vee z\quad\text{for each }k\in\Lambda\upw j\,.
 \]
In particular, $a_i^k\leq y\vee z$, thus $a_0^k\leq a_0^i\vee y\vee z=a_0^i$, for each $k\in\Lambda\upw j$. This contradicts the assumption that~$L$ has no unit.

The other possibility is $\eps_i(x)\sim_{\eps_i(z)}\eps(y)$. In such a case, $\eps_i(x)\wedge\eps_i(z)=0$, thus, \emph{a fortiori}, $\eps_i(0)=0$, that is, $a_i^k=0$ for all large enough~$k\in\Lambda$. As~$L$ has no unit, this is impossible.
\end{proof}

\begin{corollary}\label{C:BanTailId1}
Every \scml\ with a countable cofinal subset has a \Ban\ trace. Hence it embeds, as a neutral ideal and within the same quasivariety, into some \cml.
\end{corollary}

\begin{proof}
Let~$L$ be a \scml\ with an increasing cofinal sequence $\famm{e_n}{n<\omega}$. We may assume that $e_0=0$. Pick~$a_n\in L$ such that $e_n\oplus a_n=e_{n+1}$, for each $n<\omega$, and set $a_m^n:=\bigoplus\famm{a_i}{m\leq i<n}$, for all non-negative integers $m\leq n$. It is straightforward to verify that the family $\famm{a_m^n}{m\leq n<\omega}$ is a \Ban\ trace in~$L$. The second part of the statement of Corollary~\ref{C:BanTailId1} follows from Theorem~\ref{T:BanTailId}.
\end{proof}

The following definition gives an analogue, for lattices without unit, of \Ban\ functions.

\begin{definition}\label{D:BanMeas}
Let $X$ be a subset in a lattice~$L$ with zero. A \emph{$L$-valued \Ban\ measure on~$X$} is a map $\ominus\colon X^{[2]}\to L$, $(x,y)\mapsto y\ominus x$, isotone in~$y$ and antitone in~$x$, such that $y=x\oplus(y\ominus x)$ for all $x\leq y$ in~$X$.
\end{definition}

Our subsequent paper~\cite{BanCoord2} will make a heavy use of \Ban\ measures.

\begin{corollary}\label{C:BanTailId2}
Every countable \scml\ $L$ has a \Ban\ measure on~$L$.
\end{corollary}

\begin{proof}
By Corollary~\ref{C:BanTailId1}, $L$ embeds, as an ideal, into a \cml~$\tilde{L}$. Furthermore, the lattice~$\tilde{L}$ constructed in the proof of Theorem~\ref{T:BanTailId} is countable as well ($\Lambda=\omega$ is countable). By Theorem~\ref{T:BaCML}, there exists a \Ban\ function~$f$ on~$\tilde{L}$. The map $L^{[2]}\to L$, $(x,y)\mapsto y\ominus x:=y\wedge f(x)$ is obviously isotone in~$y$ and antitone in~$x$. Furthermore, it follows from the modularity of~$L$ that $y=x\oplus(y\ominus x)$ for all $x\leq y$ in~$L$. Therefore, $\ominus$ is as required.
\end{proof}

For von~Neumann regular rings the result of Corollary~\ref{C:BanRegRings} is apparently stronger.

\section{\Ban\ traces and coordinatizability}\label{S:Ban22Coord}

Coordinatizability provides another large class of lattices admitting a \Ban\ trace. 

\begin{proposition}\label{P:CoordBT}
Every coordinatizable \scml\ has a normal \Ban\ trace.
\end{proposition}

\begin{proof}
Let~$R$ be a von~Neumann regular ring, and set $\Lambda:=\Idemp R$, endowed with its ordering~$\utr$ (cf. Section~\ref{S:Basic}). Set $A_i^j:=(j-i)R$, for all $i\utr j$ in~$\Lambda$. It follows from the proof of Lemma~2 in Faith and Utumi~\cite{FaUt63} that $R$ is the directed union of its corner rings~$eRe$, where $e\in\Idemp R$. Hence, $(\Lambda,\utr)$ is upward directed and $\setm{A_0^i}{i\in\Lambda}$ is cofinal in~$\LL(R)$. It is straightforward to verify that $A_i^k=A_i^j\oplus A_j^k$ for all $i\leq j\leq k$ in~$\Lambda$. Furthermore, for $i,j\in\Lambda$ with $i\utr j$, if $A_0^i=A_0^j$, that is, $iR=jR$, then $j=ij=i$. Therefore, $\famm{A_i^j}{i\leq j\text{ in }\Lambda}$ is a normal \Ban\ trace.
\end{proof}

See also the comments following the statement of Problem~\ref{Pb:EmbId}, Section~\ref{S:Pbs}.

The following definition is taken from \cite{Herr}.

\begin{definition}\label{D:StrUnCoord}
A coordinatizable lattice~$L$ is \emph{uniquely rigidly coordinatizable}, if for all von~Neumann regular rings~$R$ and~$S$ coordinatizing~$L$, every isomorphism from~$\LL(R)$ onto~$\LL(S)$ has the form~$\LL(f)$, for a unique isomorphism $f\colon R\to S$.
\end{definition}

Hence the von~Neumann regular ring coordinatizing a uniquely rigidly coordinatizable lattice is unique up to unique isomorphism.

\begin{lemma}\label{L:1DimLiftURC}
Let $K$ be a uniquely rigidly coordinatizable principal ideal in a coordinatizable lattice~$L$, let~$R$ and~$S$ be von~Neumann regular rings with isomorphisms $\eps\colon K\to\LL(R)$ and $\eta\colon L\to\LL(S)$, and let~$e$ be an idempotent element of~$S$ such that $\eta(1_K)=eS$. Then there exists a unique ring homomorphism $f\colon R\to S$ with range~$eSe$ such that $\eta\res_K=\LL(f)\circ\eps$.
\end{lemma}

\begin{proof}
By \cite[Lemma~10.2]{Jons62}, there are mutually inverse isomorphisms
 \begin{align*}
 \alpha&\colon\LL(eSe)\to\LL(S)\dnw eS\,,\qquad J\mapsto JS\,,\\
 \beta&\colon\LL(S)\dnw eS\to\LL(eSe)\,, \qquad J\mapsto J\cap eSe\,. 
 \end{align*}
Denote by $u\colon eSe\into S$ the inclusion map and by~$\eta'$ the restriction of~$\eta$ from~$K=L\dnw 1_K$ onto $\eta(L)\dnw\eta(1_K)=\LL(S)\dnw eS$. 
We consider the following sequence of lattice embeddings:
 \[
 \xymatrix{
 \LL(R)\ar[r]^{\eps^{-1}}_{\cong}&K\ar[r]^(.3){\eta'}_(.3){\cong}&
 \LL(S)\dnw eS\ar[r]^{\beta}_{\cong}&\LL(eSe)\ar@{^(->}[r]^{\LL(u)}&\LL(S)
 }
 \]
In particular, $\beta\circ\eta'\circ\eps^{-1}\colon\LL(R)\to\LL(eSe)$ is an isomorphism, so both~$R$ and~$eSe$ coordinatize~$K$, and so, by assumption, there exists a unique isomorphism $g\colon R\onto\nobreak eSe$ such that $\LL(g)=\beta\circ\eta'\circ\eps^{-1}$. As any~$g$ satisfying $\LL(g)=\beta\circ\eta'\circ\eps^{-1}$ is necessarily one-to-one, it follows that there exists a unique surjective homomorphism $g\colon R\onto eSe$ such that $\LL(g)=\beta\circ\eta'\circ\eps^{-1}$.

As~$\beta^{-1}=\alpha$ is the restriction of~$\LL(u)$ from~$\LL(eSe)$ onto~$\LL(S)\dnw eS$, we get
 \begin{equation}\label{Eq:L(u)beetetK}
 \LL(u)\circ\beta\circ\eta'=\eta\res_K\,.
 \end{equation}
Now a ring homomorphism $f\colon R\to S$ with range $eSe$ has the form $u\circ h$, for some surjective ring homomorphism $h\colon R\onto eSe$. Then $\eta\res_K=\LL(f)\circ\eps$ if{f} $\LL(f)=(\eta\res_K)\circ\nobreak\eps^{-1}$, if{f} (using~\eqref{Eq:L(u)beetetK} together with $\LL(f)=\LL(u)\circ\LL(h)$)
$\LL(u)\circ\LL(h)=\LL(u)\circ\beta\circ\eta'\circ\eps^{-1}$, if{f} (as~$\LL(u)$ is one-to-one) $\LL(h)=\beta\circ\eta'\circ\eps^{-1}$, that is, $h=g$, which is equivalent to $f=u\circ g$.
\end{proof}

Observe that any~$f$ satisfying the condition in Lemma~\ref{L:1DimLiftURC} is necessarily an embedding from~$R$ into~$S$, so it defines by restriction an isomorphism from~$R$ onto~$eSe$. Hence the given condition on~$f$ is equivalent to the conjunction of the two following statements:
 \begin{itemize}
 \item $f$ is an embedding from~$R$ into~$S$ with range~$eSe$,
 
 \item the equality $f(x)S=(\eta\circ\eps^{-1})(xR)$ holds for each $x\in R$. 
 \end{itemize}
 
Now a variant of the argument of \cite[Theorem~10.3]{Jons62} gives the following.

\begin{proposition}\label{P:Bantr2Coord}
Let $L$ be a \scml\ with a \Ban\ trace $\famm{a_i^j}{i\leq j\text{ in }\Lambda}$ such that~$L\dnw a_0^i$ is uniquely rigidly coordinatizable for each $i\in\Lambda$. Then~$L$ is coordinatizable.
\end{proposition}

\begin{proof}
For each $i\in\Lambda$, we fix a von~Neumann regular ring~$R_i$ and an isomorphism\linebreak $\eps_i\colon L\dnw a_0^i\to\LL(R_i)$, and we denote by~$1_i$ the unit of the ring~$R_i$. For all $i\leq j$ in~$\Lambda$, it follows from the relations $R_j=\eps_j(a_0^j)=\eps_j(a_0^i)\oplus\eps_j(a_i^j)$ and Lemma~\ref{L:DecompId} that there exists a unique element $e_i^j\in\eps_j(a_0^i)$ such that $1_j-e_i^j\in\eps_j(a_i^j)$, and then
 \begin{equation}\label{Eq:Pptieseij}
 e_i^j\in\Idemp(R_j)\,,\quad\eps_j(a_0^i)=e_i^jR_j\,,\quad\text{ and }
 \eps_j(a_i^j)=(1_j-e_i^j)R_j\,.
 \end{equation}
By Lemma~\ref{L:1DimLiftURC}, there exists a unique ring embedding $f_i^j\colon R_i\into R_j$ with range $e_i^jR_je_i^j$ such that
 \begin{equation}\label{Eq:DiagrDeffij}
 \LL(f_i^j)\circ\eps_i=\eps_j\res_{L\dnw a_0^i}\,.
 \end{equation}
In particular, $f_i^j(1_i)=e_i^j$. Trivially, $f_i^i=\id_{R_i}$.

\begin{sclaim}
The equality $f_j^k(e_i^j)=e_i^k$ holds, for all $i\leq j\leq k$ in~$\Lambda$.
\end{sclaim}

\begin{scproof}
We compute
 \begin{align}
 f_j^k(e_i^j)&\in\LL(f_j^k)(e_i^jR_j)\notag\\
 &=(\LL(f_j^k)\circ\eps_j)(a_0^i)&&(\text{use~\eqref{Eq:Pptieseij}})\notag\\
 &=\eps_k(a_0^i)&&(\text{use~\eqref{Eq:DiagrDeffij}}).\label{Eq:fjkeijwhnd}
 \end{align}
Observe further that $\bigl(1_k-f_j^k(1_j)\bigr)R_k=(1_k-e_j^k)R_k=\eps_k(a_j^k)$ while
 \[
 \bigl(f_j^k(1_j)-f_j^k(e_i^j)\bigr)R_k=\LL(f_j^k)\bigl((1_j-e_i^j)R_j\bigr)=(\LL(f_j^k)\circ\eps_j)(a_i^j)=\eps_k(a_i^j)\,.
 \]
Hence,
 \begin{align}
 1_k-f_j^k(e_i^j)&=\bigl(1_k-f_j^k(1_j)\bigr)+
 \bigl(f_j^k(1_j)-f_j^k(e_i^j)\bigr)\notag\\
 &\in\eps_k(a_j^k)\oplus\eps_k(a_i^j)\notag\\
 &=\eps_k(a_i^k)\,.\label{Eq:1-fjkeijwhnd}
 \end{align}
It follows from~\eqref{Eq:fjkeijwhnd} that $f_j^k(e_i^j)\in\eps_k(a_0^i)$ while it follows from~\eqref{Eq:1-fjkeijwhnd} that\linebreak $1_k-f_j^k(e_i^j)\in\eps_k(a_i^k)$. The conclusion follows from the definition of~$e_i^k$.
\end{scproof}

Let $i\leq j\leq k$ in~$\Lambda$. It follows from the claim above that
 \[
 e_i^k\cdot e_j^k=f_j^k(e_i^j)\cdot f_j^k(1_j)=f_j^k(e_i^j\cdot 1_j)=
 f_j^k(e_i^j)=e_i^k\,,
 \]
and, similarly, $e_j^k\cdot e_i^k=e_i^k$. Hence $e_i^k\utr e_j^k$, and so
 \begin{align*}
 \im(f_j^k\circ f_i^j)&=f_j^k(e_i^jR_je_i^j)
 &&(\text{because }\im f_i^j=e_i^jR_je_i^j)\\
 &=f_j^k(e_i^j)\bigl(e_j^kR_ke_j^k\bigr)f_j^k(e_i^j)
 &&(\text{because }\im f_j^k=e_j^kR_ke_j^k)\\
 &=e_i^ke_j^kR_ke_j^ke_i^k&&(\text{by the claim above})\\
 &=e_i^kR_ke_i^k\,. 
 \end{align*}
Now for each $x\in R_i$, it follows from~\eqref{Eq:DiagrDeffij} that $f_i^j(x)R_j=(\eps_j\circ\eps_i^{-1})(xR_i)$, while, setting $y:=f_i^j(x)$, we get
$f_j^k(y)R_k=(\eps_k\circ\eps_j^{-1})(yR_j)$, so
 \[
 (f_j^k\circ f_i^j)(x)R_k=f_j^k(y)R_k=
 (\eps_k\circ\eps_j^{-1}\circ\eps_j\circ\eps_i^{-1})(xR_i)=
 (\eps_k\circ\eps_i^{-1})(xR_i)\,.
 \]
Therefore, by the uniqueness of the property defining~$f_i^k$, we obtain that the equality $f_i^k=f_j^k\circ f_i^j$ holds.

It follows that we can form the direct limit
 \[
 \famm{R,f_i}{i\in\Lambda}=\varinjlim\famm{R_i,f_i^j}{i\leq j\text{ in }\Lambda}\,.
 \]
As~$R$ is a direct limit of von~Neumann regular rings, it is a von~Neumann regular ring. As the functor~$\LL$ preserves direct limits, we obtain that
 \[
 L=\varinjlim_{i\in\Lambda}(L\dnw a_0^i)\cong\varinjlim_{i\in\Lambda}\LL(R_i)
 \cong\LL(R)\,,
 \]
and so~$L$ is coordinatizable.
\end{proof}

\begin{remark}\label{Rk:Bantr2Coord}
The example, presented at the bottom of Page~301 in \cite{Jons62}, of the lattice of all finite-dimensional subspaces of a vector space of countable infinite dimension, shows that the conclusion of Proposition~\ref{P:Bantr2Coord} cannot be strengthened to saying that~$L$ is uniquely coordinatizable.
\end{remark}

\begin{theorem}\label{T:CharCoord4Fr}
Let~$L$ be a \scml\ that admits a large $4$-frame, or which is Arguesian and that admits a large $3$-frame. Then the following are equivalent:
\begin{enumerate}
\item $L$ is coordinatizable;

\item $L$ has a normal \Ban\ trace;

\item $L$ has a \Ban\ trace.
\end{enumerate}
\end{theorem}

\begin{proof}
The direction (i)$\Rightarrow$(ii) follows from Proposition~\ref{P:CoordBT}, while (ii)$\Rightarrow$(iii) is trivial.

Now let $L$ be a \scml\ with a large $n$-frame $(\famm{a_s}{0\leq s<n},\famm{c_s}{1\leq s<n})$, where $n\geq 4$, or only $n\geq 3$ in case $L$ is Arguesian; set $a:=\bigvee_{s<n}a_s$. If $L$ has a \Ban\ trace $\famm{e_i^j}{i\leq j\text{ in }\Lambda}$, then we may assume, replacing~$\Lambda$ by~$\Lambda\upw i_0$ for an index~$i_0$ such that $a\leq e_0^{i_0}$, that the inequality $a\leq e_0^i$ holds for each $i\in\Lambda$. As the element~$a$ is large in~$L$, it follows easily from \cite[Lemma~1.4]{Jons60} that~$a$ is large in each~$L\dnw e_0^i$ as well. 

Now it is observed in \cite[Theorem~10.4]{Jons62} that every \cml\ that admits a large $4$-frame, or which is Arguesian and that admits a large $3$-frame, is uniquely coordinatizable; the conclusion is strengthened to ``uniquely rigidly coordinatizable'' in~\cite[Corollary~4.12]{Micol}, see also~\cite[Theorem~18]{Herr}. In particular, all the lattices $L\dnw e_0^i$, for $i\in\Lambda$, are uniquely rigidly coordinatizable. Therefore, by Proposition~\ref{P:Bantr2Coord}, $L$ is coordinatizable.
\end{proof}

We shall prove in \cite{BanCoord2} that there exists a non-coordinatizable \scml~$L$ with a large $4$-frame. Hence~$L$ does not have a \Ban\ trace as well. The construction of~$L$ requires techniques far beyond those involved in the present paper.

\section{Problems}\label{S:Pbs}

By Fuchs and Halperin~\cite{FuHa64}, every von~Neumann regular ring~$R$ can be embedded as a two-sided ideal into some \emph{unital} von~Neumann regular ring~$S$. Consequently, $\LL(R)$ embeds as a neutral ideal into~$\LL(S)$. This gives a proof, that uses neither Theorem~\ref{T:BanTailId} nor Proposition~\ref{P:CoordBT}, that every coordinatizable \scml\ embeds as a neutral ideal into some coordinatizable \cml. We do not know the general answer in the non-coordinatizable case:

\begin{problem}\label{Pb:EmbId}
Does every \scml\ embed as a (neutral) ideal into some \cml?
\end{problem}

It is proved in Theorem~\ref{T:BaCML} that every countable \cml\ has a Boolean \Ban\ function. The range of such a \Ban\ function is easily seen to be a \emph{maximal} Boolean sublattice of~$L$.

\begin{problem}\label{Pb:MaxBool}
Is every maximal Boolean sublattice of a countable \cml~$L$ the range of some \Ban\ function on~$L$? Are any two such Boolean sublattices isomorphic?
\end{problem}

Finally, we should mention that while the present paper is devoted to \emph{modular} lattices, the notion of a \Ban\ function is also well-defined for non-modular lattices.

\begin{problem}\label{Pb:BanFctRelCpl}
Does every countable, bounded, relatively complemented lattice have a \Ban\ function?
\end{problem}

Observe that Example~\ref{Ex:NonBanCL} gives a finite complemented lattice without a \Ban\ function. Also observe that the existence of a \Ban\ function on a bounded lattice~$L$ does not imply in general that~$L$ is relatively complemented, which suggests that Problem~\ref{Pb:BanFctRelCpl} may not be the ``right'' question.

\section*{Acknowledgment}
I thank deeply Christina Brzuska, Ken Goodearl, and Pavel R\r{u}\v{z}i\v{c}ka for their enlightening comments, that lead to many improvements of the paper, and in particular to a much more transparent proof of Theorem~\ref{T:BaCML}. I also thank the referees for their careful reading of the paper.


\begin{thebibliography}{99}
\bibitem{Bana}
B. Banaschewski,
\emph{Totalgeordnete Moduln} (German),
Arch. Math. \textbf{7} (1957), 430--440.

\bibitem{Birk94}
G. Birkhoff,
``Lattice Theory''. Corrected reprint of the 1967 third edition. American Mathematical Society Colloquium Publications~\textbf{25}. American Mathematical Society, Providence, R.I., 1979. vi+418~p. ISBN: 0-8218-1025-1\,.

\bibitem{Dobb83}
H. Dobbertin,
\emph{Refinement monoids, Vaught monoids, and Boolean algebras}, 
Math. Ann.~\textbf{265}, no.~4 (1983), 473--487.

\bibitem{FaUt63}
C. Faith and Y. Utumi,
\emph{On a new proof of Litoff's theorem},
Acta Math. Acad. Sci. Hungar.~\textbf{14} (1963), 369--371.

\bibitem{FrHa54}
K.\,D. Fryer and I. Halperin,
\emph{Coordinates in geometry},
Trans. Roy. Soc. Canada. Sect. III. (3)~\textbf{48} (1954), 11--26.

\bibitem{FuHa64}
L. Fuchs and I. Halperin,
\emph{On the imbedding of a regular ring in a regular ring with identity}, Fund. Math.~\textbf{54} (1964), 285--290.

\bibitem{Good91}
K.\,R. Goodearl,
``Von Neumann Regular Rings'',
Second edition.  Robert E. Krieger Publishing Co., Inc., Malabar, FL, 1991.
xviii+412~p. ISBN: 0-89464-632-X\,.

\bibitem{GLT2}
G. Gr\"atzer,
``General Lattice Theory. Second edition'', new appendices by the
author with B.\,A. Davey, R. Freese, B. Ganter,
M. Greferath, P. Jipsen, H.\,A. Priestley, H. Rose, E.\,T. Schmidt,
S.\,E. Schmidt, F. Wehrung, and R. Wille. Birkh\"auser Verlag, Basel,
1998. xx+663~p. ISBN: 3-7643-5239-6 (Basel), 0-8176-5239-6 (Boston).

\bibitem{Herr}
C. Herrmann,
\emph{Generators for complemented modular lattices and  
the von~Neu\-mann-J\'{o}nsson Coordinatization Theorems}, Algebra Universalis, to appear.

\bibitem{Jons60}
B. J\'onsson,
\emph{Representations of complemented modular lattices},
Trans. Amer. Math. Soc.~\textbf{60} (1960), 64--94.

\bibitem{Jons62}
B. J\'onsson,
\emph{Representations of relatively complemented modular lattices},
Trans. Amer. Math. Soc.~\textbf{103} (1962), 272--303.

\bibitem{Maed58}
F. Maeda,
``Kontinuierliche Geometrien''. (German) 
Die Grundlehren der mathematischen Wissenschaften in Einzeldarstellungen mit besonderer Ber\"ucksichtigung der Anwendungsgebiete, Bd.~\textbf{95}. Springer-Verlag, Berlin-G\"ottingen-Heidelberg, 1958. x+244~p.

\bibitem{Micol}
F.~Micol,
``On representability of $\ast$-regular rings and modular ortholattices'',
PhD thesis, TU Darmstadt, January 2003. Available online at\newline \texttt{http://elib.tu-darmstadt.de/diss/000303/diss.pdf}\,.

\bibitem{HBA3}
J.\,D. Monk (ed.),
``Handbook of Boolean algebras'', Volume~\textbf{3}. Ed. with the cooperation of R. Bonnet. (English) North-Holland, Amsterdam - New York - Oxford - Tokyo, 1989. xix+1367~p. ISBN: 0 444 87153 5\,.

\bibitem{Pier}
R.\,S. Pierce,
\emph{Countable Boolean algebras}, Chapter~21 in~\cite{HBA3}, 775--876.

\bibitem{SaSo}
M. Saarim\"aki and P. Sorjonen,
\emph{On Banaschewski functions in lattices}, Algebra Universalis \textbf{28}, no.~1 (1991), 103--118.

\bibitem{WDim}
F. Wehrung, \emph{The dimension monoid of a lattice}, Algebra Universalis~\textbf{40}, no.~3 (1998), 247--411.

\bibitem{CXCoord}
F. Wehrung,
\emph{Von Neumann coordinatization is not first-order}, J. Math. Log.~\textbf{6}, no.~1 (2006), 1--24.

\bibitem{BanCoord2}
F. Wehrung,
\emph{A non-coordinatizable \scml\ with a large J\'onsson four-frame}, preprint.

\end{thebibliography}
\end{document}